\documentclass[a4paper,12pt]{article}
\usepackage{amssymb}
\usepackage{amsfonts}
\usepackage{amsmath}
\usepackage{amsthm}
\newtheorem{defin}{Definition}[section]
\newtheorem{prop}[defin]{Proposition}
\newtheorem{cor}[defin]{Corollary}
\newtheorem{lemma}[defin]{Lemma}
\newcommand{\recht}{\rightarrow}
\newcommand{\hilh}{\mathcal{H}}
\newcommand{\rondm}{\mathcal{M}}
\newcommand{\gewn}{\mathfrak{N}}
\newcommand{\gewm}{\mathfrak{M}}
\newcommand{\reeel}{\mathbb{R}}
\newcommand{\complex}{\mathbb{C}}
\newcommand{\ffi}{\varphi}
\newcommand{\ffid}{\varphi_\delta}
\newcommand{\et}{\mbox{\tiny $\frac{1}{2}$}}
\newcommand{\iet}{\mbox{\tiny $\frac{i}{2}$}}
\newcommand{\rondz}{\mathcal{Z}}

\begin{document}
\title{A Radon-Nikodym theorem for von Neumann algebras}
\author{Stefaan Vaes\footnote{Research Assistant of the Fund for Scientific
Research - Flanders (Belgium)(F.W.O.)}\\
Department of Mathematics \\Katholieke Universiteit Leuven
\\Celestijnenlaan 200B\\3001 Heverlee\\Belgium\\email : Stefaan.Vaes@wis.kuleuven.ac.be
\\fax : (32)-16 32 79 98}
\date{\textbf{November 1998}}
\maketitle
\begin{abstract}
\noindent In this paper we present a generalization of the
Radon-Nikodym theorem proved by Pedersen and Takesaki in
\cite{pedersen}. Given a normal, semifinite and faithful (n.s.f.)
weight $\ffi$ on a von Neumann algebra $\rondm$ and a strictly
positive operator $\delta$, affiliated with $\rondm$ and satisfying a
certain relative invariance property with respect to the modular
automorphism group $\sigma^\ffi$ of $\ffi$, with a strictly positive
operator as the invariance factor, we construct the n.s.f. weight
$\varphi(\delta^{\et} \cdot \delta^{\et})$. All the n.s.f. weights on
$\rondm$ whose modular automorphisms commute with $\sigma^\ffi$ are
of this form, the invariance factor being affiliated with the centre of
$\rondm$. All the n.s.f. weights which are relatively invariant under
$\sigma^\ffi$ are of this form, the invariance factor being a scalar.
\end{abstract}
\section*{Introduction}
In \cite{pedersen} G.K. Pedersen and M. Takesaki gave a construction
of a normal semifinite faithful (n.s.f.) weight $\varphi(\, \cdot \,
\delta)$ on a von Neumann algebra $\mathcal{M}$, starting from a
n.s.f. weight $\varphi$ on $\mathcal{M}$ and a strictly positive
operator $\delta$ affiliated with the von Neumann algebra of elements
invariant under the modular automorphisms of $\varphi$. These weights
$\varphi(\, \cdot \,
\delta)$ are precisely all the n.s.f. weights $\psi$ on $\mathcal{M}$
which are invariant under the modular automorphisms of $\varphi$. In this
paper we will give a construction for a n.s.f. weight
$\varphi(\delta^{\mbox{\tiny $\frac{1}{2}$}} \cdot \delta^{\mbox{\tiny $\frac{1}{2}$}})$ in the case where $\delta$
satisfies the weaker hypothesis $\sigma_s^\varphi(\delta^{it}) =
\lambda^{ist} \delta^{it}$ for all $s,t \in \mathbb{R}$ and for a
given strictly positive operator $\lambda$, affiliated with
$\mathcal{M}$ and strongly commuting with $\delta$. This way we obtain
precisely all n.s.f. weights $\psi$ on $\mathcal{M}$ for which $[D\psi
: D \varphi]_t = \lambda^{\mbox{\tiny $\frac{1}{2}$} it^2} \delta^{it}$.
The operators $\lambda$ and $\delta$ are uniquely determined by
$\psi$.
When $\psi$ is a n.s.f. weight on $\rondm$ we prove that
$\sigma^\psi$ and $\sigma^\ffi$ commute if and only if there exist
strictly positive operators $\lambda$ and $\delta$ affiliated with
the centre of $\rondm$ and $\rondm$ respectively, such that
$\sigma_s^\ffi(\delta^{it}) = \lambda^{ist} \delta^{it}$ for all $s,t
\in \reeel$ and such that $\psi = \varphi(\delta^{\et} \cdot
\delta^{\et})$. When $\psi$ is a n.s.f. weight on $\rondm$ and
$\lambda \in \reeel_0^+$ we prove that $\psi \circ \sigma_t^\ffi =
\lambda^{-t} \psi$ for all $t \in \reeel$ if and only if $\ffi \circ
\sigma_t^\psi = \lambda^t \psi$ for all $t \in \reeel$ if and only if
there exists a strictly positive operator $\delta$ affiliated with $\rondm$
such that $\sigma_s^\ffi(\delta^{it}) = \lambda^{ist} \delta^{it}$
for all $s,t \in \reeel$ and such that $\psi =\varphi(\delta^{\et} \cdot
\delta^{\et})$.

One important application of the Radon-Nikodym theorem of Pedersen and
Takesaki arose in the theory of locally compact quantum groups. In
\cite{masuda} the theorem is used to obtain the modular element as the
Radon-Nikodym derivative of the left and the right Haar weight. Very
recently the new and very simple definition of locally compact
quantum groups we found in cooperation with J.~Kustermans (see
\cite{kustvaes}) implies that in general the right Haar weight is
only relatively invariant under the modular automorphisms of the left
Haar weight. In order
still to be able to obtain the modular element we need the more
general Radon-Nikodym theorem of this paper. This is a very important
application. Further, the possibility to obtain a Radon-Nikodym
derivative from the sole assumption that $\sigma^\ffi$ and
$\sigma^\psi$ commute is new and could give rise to several
applications in von Neumann algebra theory. It is also very important
to notice that the most powerful tool to prove the equality of two
n.s.f. weights, namely showing that the Radon-Nikodym derivative is
trivial, can now be applied in much more situations.

Let us first fix some notations. In paragraphs 1 to 4 we will always
assume that $\mathcal{M}$ is a von Neumann algebra, that $\varphi$ is
a n.s.f. weight on $\mathcal{M}$ and that $\lambda$ and $\delta$ are two
strictly positive, strongly commuting operators affiliated with
$\mathcal{M}$. We suppose $\mathcal{M}$ acts on the GNS-space
$\mathcal{H}$ of $\varphi$. We denote by $J$ and $\Delta$ the modular
operators of $\ffi$ and by $(\sigma_t)$ the modular automorphisms. As
usual we put $\gewn = \{a \in \rondm \mid \ffi(a^*a) < \infty \}$ and
$\gewm=\gewn^* \gewn$. We denote by $\Lambda : \gewn \recht \hilh$
the map appearing in the GNS-construction of $\ffi$ such that
$\langle \Lambda(a), \Lambda(b) \rangle = \ffi(b^*a)$. We remark that
the map $\Lambda$ is weak operator -- weak closed and refer to
\cite{stratila}, Chapter 10 and \cite{stratila2}, Chapter I,
for more details about n.s.f. weights.
We assume the following relative invariance~:
\begin{equation*}
\sigma_t(\delta^{is}) = \lambda^{ist}\delta^{is} \quad\text{for
all}\quad s,t \in \reeel.
\end{equation*}
Remark that in case $\lambda = 1$ we arrive at the premises for the
construction of Pedersen and Takesaki. Because we will regularly use
analytic continuations, we introduce the notation $S(z)$ for the
closed strip of complex numbers with real part between $0$ and
$\text{Re}(z)$.

Starting from all these assumptions we will construct a n.s.f. weight
$\ffid$ on $\rondm$ in the first paragraph. Then we will compute the
modular operators and automorphisms of $\ffid$ and prove an explicit
formula that justifies the notation $\ffid=\ffi(\delta^{\et} \, \cdot
\, \delta^{\et})$. In the fourth paragraph we compute the Connes
cocycle $[D \ffid : D \ffi ]$, which will enable us to prove in the last paragraph
the three Radon-Nikodym type theorems mentioned above.
\section{The construction of the weight $\ffid$}
\begin{defin}
For each $n \in \mathbb{N}_0$ we define an element $e_n \in \rondm$ by
\begin{align*}
\alpha_n &= \frac{2n^2}{\Gamma(\frac{1}{2})\Gamma(\frac{1}{4})},
\\
e_n &= \alpha_n \int_{-\infty}^{+\infty}\int_{-\infty}^{+\infty}
\exp(-n^2 x^2-n^4 y^4) \lambda^{ix} \delta^{iy} \; dx \; dy \in \rondm.
\end{align*}
\end{defin}
The integral makes sense in the strong* topology. We remark that
automotically $\lambda$ satisfies $\sigma_t(\lambda^{is})=\lambda^{is}$ for all
$s,t \in \reeel$, which can be proven very easily. We also easily
obtain the following lemma, using the 'analytic extension techniques'
of \cite{stratila}, Chapter 9.
\renewcommand{\theenumi}{\roman{enumi}}
\renewcommand{\labelenumi}{\theenumi)}
\begin{lemma}
\begin{enumerate}
\item The elements $e_n \in \rondm$ are analytic w.r.t. $\sigma$.
For all $x,y,z \in \complex$ the operator $\delta^x \lambda^y
\sigma_z(e_n)$ is bounded, with domain $\hilh$, analytic w.r.t.
$\sigma$
and satisfies $\sigma_t(\delta^x \lambda^y \label{en1}
\sigma_z(e_n)) = \delta^x \lambda^{y+tx}
\sigma_{t+z}(e_n)$ for all $t \in \complex$.
\item For all $z \in \complex$ we have $\sigma_z(e_n) \recht 1$
strong* and bounded. \label{en2}
\item The function $(x,y,z) \mapsto \delta^x \lambda^y
\sigma_z(e_n)$ is analytic from $\complex^3$ to $\rondm$. \label{en3}
\item The elements $e_n$ are selfadjoint.
\end{enumerate}
\end{lemma}
Inspired by the work of \cite{kustermans1} we give the following
definition :
\begin{defin}
Define a subset $\gewn_0$ of $\rondm$ by
\begin{align*}
&\gewn_0 = \{a \in \rondm \mid a \delta^{\mbox{\tiny $\frac{1}{2}$}} \quad\text{is bounded and}\quad
\overline{a\delta^{\mbox{\tiny $\frac{1}{2}$}}} \in \gewn\} \\
\intertext{and a map}
&\Gamma : \gewn_0 \recht \hilh : a \mapsto
\Lambda(\overline{a\delta^{\mbox{\tiny $\frac{1}{2}$}}}),
\end{align*}
where $\overline{a \delta^{\et}}$ denotes the closure of $a
\delta^{\et}$.
\end{defin}
Remark that $\Gamma$ is injective and $\gewn_0$ is a left ideal in
$\rondm$. So $\Gamma(\gewn_0 \cap \gewn_0^*)$ becomes an involutive algebra by
defining
\begin{align*}
\Gamma(a)\Gamma(b) &= \Gamma(ab) \\
\Gamma(a)^\# &= \Gamma(a^*).
\end{align*}
\begin{prop}
When we endow the involutive algebra $\Gamma(\gewn_0 \cap \gewn_0^*)$
with the scalar product of $\hilh$, it becomes
a left Hilbert algebra. The generated von Neumann algebra is $\rondm$.
\end{prop}
\begin{proof}
If $a,b \in \gewn_0 \cap \gewn_0^*$ we have
\begin{equation*}
\Gamma(ab) = \Lambda(a \overline{b \delta^{\mbox{\tiny $\frac{1}{2}$}}}) = a \Gamma(b)
\end{equation*}
so that $\Gamma(b) \mapsto \Gamma(ab)$ is bounded. For $a,b,c \in \gewn_0
\cap \gewn_0^*$ we have
\begin{equation*}
\langle \Gamma(a) \Gamma(b),\Gamma(c) \rangle = \varphi((\overline{c \delta^{\mbox{\tiny
$\frac{1}{2}$}}})^*
a (\overline{b \delta^{\mbox{\tiny $\frac{1}{2}$}}})) = \langle \Gamma(b),
\Lambda(a^* (\overline{c \delta^{\mbox{\tiny $\frac{1}{2}$}}})) \rangle = \langle \Gamma(b),
\Gamma(a)^\# \Gamma(c) \rangle.
\end{equation*}

If $a \in \gewn \cap \gewn^*$ one can easily verify that $e_n a
(\delta^{-\mbox{\tiny $\frac{1}{2}$}}e_n) \in \gewn_0 \cap \gewn_0^*$. Moreover
\begin{equation*}
\Gamma(e_n a (\delta^{-\mbox{\tiny $\frac{1}{2}$}}e_n)) = \Lambda(e_nae_n) =
J(\sigma_{\iet}(e_n))^*Je_n \Lambda(a) \recht
\Lambda(a).
\end{equation*}
So $\Gamma(\gewn_0 \cap \gewn_0^*)$ is dense in $\hilh$. But also $e_n a e_n
\in \gewn_0 \cap \gewn_0^*$, and this converges strongly to $a$.
Therefore $\gewn_0 \cap \gewn_0^*$ is strongly dense in $\rondm$ and thus
$(\Gamma(\gewn_0 \cap \gewn_0^*))^2$ is dense in $\hilh$.

We claim that for all $n \in \mathbb{N}_0$ and all $b$ and $b'$ in the
Tomita algebra of $\ffi$, the element $\Lambda(e_n b b' e_n)$ belongs
to the domain of the adjoint of the mapping $\Gamma(a) \mapsto
\Gamma(a)^\#$. This will imply the closedness of that mapping, and so
this will end the proof. To prove the claim, choose $a \in \gewn_0 \cap \gewn_0^*$. Define
the element $x \in \gewn$ by
\begin{equation*}
x:= (\delta^{\mbox{\tiny $\frac{1}{2}$}} \sigma_{-i}(e_n))
\sigma_{-i}({b'}^*b^*(\delta^{-\et}e_n)).
\end{equation*}
Then we can make the following calculation~:
\begin{align*}
\langle \Lambda(x),\Gamma(a) \rangle &=\varphi((\overline{a \delta^{\et}})^*
(\delta^{\et} \sigma_{-i}(e_n))
\sigma_{-i}({b'}^*) \sigma_{-i}(b^*(\delta^{-\et}e_n))) \\
&=\varphi(b^*(\delta^{-\et}e_n) (\overline{a \delta^{\et}})^*
(\delta^{\et} \sigma_{-i}(e_n)) \sigma_{-i}({b'}^*)) \\
&=\varphi(b^*e_na^*(\delta^{\et}\sigma_{-i}(e_n))\sigma_{-i}({b'}^*)) \\
&= \varphi(b^* e_n (\overline{a^* \delta^{\et}}) \sigma_{-i}(e_n{b'}^*))
\\ &=\varphi(e_n {b'}^*b^*e_n (\overline{a^* \delta^{\et}})) \\
&=\langle \Gamma(a)^\#,\Lambda(e_nbb'e_n) \rangle.
\end{align*}
This proves our claim.
\end{proof}
\begin{defin} \label{def15}
We define $\ffid$ as the weight associated to the left Hilbert
algebra\linebreak
$\Gamma(\gewn_0 \cap \gewn_0^*)$. This is a n.s.f. weight on
$\rondm$.
\end{defin}
We denote by $\gewn', \gewm'$ and $\Lambda' : \gewn' \recht \hilh$
the evident objects associated to $\ffid$. We denote by
$(\sigma'_t)$ the modular automorphisms of $\ffid$. We remark that
$\gewn_0 \subset \gewn'$ and $\Lambda'(a) = \Gamma(a)$ for all $a \in
\gewn_0$.

Up to now the operator $\lambda$ did not appear in our formulas. We
only need the relative invariance property of $\delta$ to construct
the analytic elements $e_n$, which cut down $\delta$ properly.
Further on $\lambda$ will of course appear when we prove properties
about $\ffid$.
\section{The modular operators of $\ffid$}
We will now calculate the modular operators and the modular
automorphisms of $\ffid$. We will give explicit formulas.
\begin{lemma}
For all $s \in \reeel$ define
\begin{equation*}
u_s = J \lambda^{\et is^2} \delta^{is} J
\lambda^{\et is^2} \delta^{is} \Delta^{is}.
\end{equation*}
Then $(u_s)$ is a strongly continuous one-parameter group of unitaries on
$\hilh$.
\end{lemma}
\begin{proof}
Straightforward, by using the facts that $J \rondm J = \rondm'$,
$J\Delta^{is} = \Delta^{is}J$, $\Delta^{is} \delta^{it} =
\lambda^{ist}\delta^{it} \Delta^{is}$ and $\Delta^{is} \lambda^{it} =
\lambda^{it} \Delta^{is}$ for all $s,t \in \reeel$.
\end{proof}
\begin{defin}
We define $\Delta'$ as the strictly positive operator on $\hilh$ such
that $u_s = {\Delta'}^{is}$ for all $s \in \reeel$.
\end{defin}
Further on, in proposition~\ref{prop24}, we will give a more explicit formula for $\Delta'$. We
first need a lemma that we will use several times.
\begin{lemma} \label{lemma8}
Let $z \in \complex$ and $n,m \in \mathbb{N}_0$. \\ If $\xi \in
\mathcal{D}({\Delta'}^z)$ then $Je_nJe_m \xi \in \mathcal{D}({\Delta'}^z)
\cap \mathcal{D}(\Delta^z)$ and
\begin{align*}
{\Delta'}^z Je_nJe_m \xi &= J \sigma_{i \bar{z}}(e_n)J \sigma_{-i
z}(e_m) \; {\Delta'}^z \xi \\
\Delta^z Je_nJe_m \xi &= J \lambda^{\et i \bar{z}^2}
\delta^{\bar{z}} \sigma_{i \bar{z}}(e_n)J \; \lambda^{\et iz^2}
\delta^{-z} \sigma_{-iz}(e_m) \; {\Delta'}^z \xi.
\end{align*}
If $\xi \in \mathcal{D}(\Delta^z)$ then $Je_nJe_m \xi \in \mathcal{D}({\Delta'}^z)
\cap \mathcal{D}(\Delta^z)$ and
\begin{align*}
{\Delta'}^z Je_nJe_m \xi &= J \lambda^{-\et i \bar{z}^2}
\delta^{-\bar{z}} \sigma_{i \bar{z}}(e_n)J \; \lambda^{-\et iz^2}
\delta^{z} \sigma_{-iz}(e_m) \; \Delta^z \xi \\
\Delta^z Je_nJe_m \xi &=J \sigma_{i \bar{z}}(e_n)J \sigma_{-i
z}(e_m) \; \Delta^z \xi.
\end{align*}
\end{lemma}
\begin{proof}
Let $\xi \in \mathcal{D}({\Delta'}^z)$. Recall the notation $S(z)$
from the end of the introduction. We define the function from $S(z)$
to $\hilh$ that maps $\alpha$ to
\begin{equation*}
J \lambda^{\et i
\bar{\alpha}^2}
\delta^{\bar{\alpha}} \sigma_{i \bar{\alpha}}(e_n)J \; \lambda^{\et i\alpha^2}
\delta^{-\alpha} \sigma_{-i\alpha}(e_m) \; {\Delta'}^\alpha \xi.
\end{equation*}
This function is continuous on $S(z)$ and analytic on its interior. In $is$ it attains the value
\begin{equation*}
J \lambda^{-\et is^2} \delta^{-is} \sigma_s(e_n)J
\; \: \lambda^{-\et is^2} \delta^{-is} \sigma_s(e_m) \; {\Delta'}^{is}
\xi \quad
= \quad J\sigma_s(e_n)J \; \sigma_s(e_m) \; \Delta^{is} \xi \quad = \quad \Delta^{is} \; Je_nJe_m
\xi.
\end{equation*}
By the results of \cite{stratila}, Chapter 9, the second statement follows. The
three remaining statements are proved analogously.
\end{proof}
\begin{prop} \label{prop24}
Let $r \in \reeel$. The operator
\begin{equation*}
J \lambda^{-\et ir^2}J \; \lambda^{-\et ir^2} \; J
\delta^{-r}J \; \delta^r \; \Delta^r
\end{equation*}
is closable and its closure equals ${\Delta'}^r$.
\end{prop}
\begin{proof}
Let $\xi \in \mathcal{D}(J
\delta^{-r}J  \; \delta^r \; \Delta^r)$. Let $n,m \in \mathbb{N}_0$. By
lemma~\ref{lemma8} we have $Je_nJe_m \xi \in \mathcal{D}({\Delta'}^r)$
and
\begin{equation*}
{\Delta'}^rJe_nJe_m \xi = J \lambda^{-\et ir^2} \sigma_{ir}(e_n)J
\lambda^{-\et ir^2} \sigma_{-ir}(e_m)
\; J
\delta^{-r}J  \delta^r \Delta^r \xi.
\end{equation*}
The operator ${\Delta'}^r$ being closed, we obtain that $\xi \in
\mathcal{D}({\Delta'}^r)$ and
\begin{equation*}
{\Delta'}^r \xi = J \lambda^{-\et ir^2}J
\lambda^{-\et ir^2} \; J
\delta^{-r}J  \delta^r \Delta^r \xi.
\end{equation*}
On the other hand let $\xi \in \mathcal{D}({\Delta'}^r)$. Let $n,m \in
\mathbb{N}_0$. By lemma~\ref{lemma8} we have that $Je_nJe_m \xi \in
\mathcal{D}(J
\delta^{-r}J  \delta^r \Delta^r)$ and ${\Delta'}^rJe_nJe_m \xi \recht
{\Delta'}^r \xi$. This implies that $\mathcal{D}(J
\delta^{-r}J  \delta^r \Delta^r)$ is a core for ${\Delta'}^r$, and
this ends our proof.
\end{proof}
Denote by $S'$ the closure of the operator $\Gamma(a) \mapsto
\Gamma(a)^\#$ on $\Gamma(\gewn_0 \cap \gewn_0^*)$. Define $J' = J
\lambda^{-i/8} J \lambda^{i/8}J$.
\begin{prop}
\begin{equation*}
S'=J' {\Delta'}^{\et}.
\end{equation*}
So, $J'$ and $\Delta'$ are the modular operators associated with
$\ffid$.
\end{prop}
\begin{proof}
Let $a \in \gewn_0 \cap \gewn_0^*$ and $n,m,k,l \in \mathbb{N}_0$. Then
$\Lambda(e_k a (\delta^{\et} e_l)) \in \mathcal{D}(\Delta^{\et})$, so
by lemma~\ref{lemma8} we have
\begin{equation*}
Je_nJe_m \Lambda(e_k a (\delta^{\et}e_l)) \in
\mathcal{D}({\Delta'}^{\et})
\end{equation*}
and
\begin{align*}
J'{\Delta'}^{\et} \; Je_nJe_m \Lambda(e_k a (\delta^{\et}e_l)) &=
\delta^{-\et}\sigma_{\iet}(e_n) \; J \lambda^{-i/4} \delta^{\et}
\sigma_{-\iet}(e_m) \; \Delta^{\et} \Lambda(e_k a (\delta^{\et}e_l)) \\
&=\delta^{-\et}\sigma_{\iet}(e_n) \; J \sigma_{-\iet}(\delta^{\et}e_m)J
\; \Lambda((\delta^{\et}e_l)a^* e_k) \\
&=\sigma_{\iet}(e_n)e_l \; \Lambda(a^*(\delta^{\et}e_m)e_k) \\
&=\sigma_{\iet}(e_n)e_l \; J \sigma_{-\iet}(e_m e_k)J \; \Gamma(a^*).
\end{align*}
The last expression converges to $\Gamma(a^*)=S' \Gamma(a)$, while
\begin{equation*}
Je_nJe_m \Lambda(e_k a (\delta^{\et}e_l)) = Je_n
\sigma_{-\iet}(e_l)Je_m e_k \Gamma(a)
\end{equation*}
converges to $\Gamma(a)$ when $n,m,k,l \recht \infty$. This implies
that $\Gamma(a) \in \mathcal{D}({\Delta'}^{\et})$ and $J' {\Delta'}^{\et}\Gamma(a) = S'\Gamma(a)$.
Thus, $S' \subset J' {\Delta'}^{\et}$.

On the other hand let $\xi \in \mathcal{D}(J \delta^{-\et}J
\delta^{\et} \Delta^{\et})$. Take a sequence $(\xi_k)$ in $\Lambda(\gewn \cap
\gewn^*)$ such that $\xi_k \recht \xi$ and $\Delta^{\et}\xi_k \recht
\Delta^{\et} \xi$. Let $n,m,k \in \mathbb{N}_0$. Then $Je_nJe_m \xi_k
\in \mathcal{D}({\Delta'}^{\et})$ and
\begin{equation*}
{\Delta'}^{\et}Je_nJ e_m \xi_k = J \lambda^{-i/8} \delta^{-\et}
\sigma_{\iet}(e_n)J \lambda^{-i/8}\delta^{\et} \sigma_{-\iet}(e_m)
\Delta^{\et} \xi_k.
\end{equation*}
If $k \recht \infty$ this converges to
\begin{equation*}
J \lambda^{-i/8} \sigma_{\iet}(e_n)J \lambda^{-i/8}
\sigma_{-\iet}(e_m) \;
J \delta^{-\et}J \delta^{\et} \Delta^{\et} \xi
= J \sigma_{\iet}(e_n)J \sigma_{-\iet}(e_m){\Delta'}^{\et} \xi.
\end{equation*}
If $n,m \recht \infty$ this converges to ${\Delta'}^{\et} \xi$.
Because $Je_nJe_m \xi_k \in \mathcal{D}(S')$ for all $n,m,k \in
\mathbb{N}$ and because of the previous proposition,
we have finally proved that $\mathcal{D}(S')$ is a core
for ${\Delta'}^{\et}$.
\end{proof}
\begin{cor} \label{Cor}
We have the formula
\begin{equation*}
\sigma'_s(x) = \lambda^{\et is^2} \delta^{is} \sigma_s(x)
\delta^{-is} \lambda^{-\et is^2}
\end{equation*}
for all $s \in \reeel$ and all $x \in \rondm$.
\end{cor}
\begin{cor}
For all $s \in \reeel$, $x,y,z \in \complex$ and $n \in \mathbb{N}_0$
we have
\begin{equation*}
\sigma_s(\lambda^x \delta^y \sigma_z(e_n)) = \sigma'_s(\lambda^x
\delta^y \sigma_z(e_n)).
\end{equation*}
\end{cor}
Remark that formulas become easier when $\lambda$ is affiliated with the
centre of $\rondm$, in particular when $\lambda$ is a positive real
number. In that case $\Delta'$ is the closure of $J \delta^{-1} J
\delta \Delta$ and $J'$ equals $\lambda^{i/4}J$, because $J x = x^* J$
for all $x$ belonging to the centre of $\rondm$. Moreover we have
$\sigma'_s(x) = \delta^{is} \sigma_s(x) \delta^{-is}$ in that case.
\section{A formula for $\ffid$}
Before we can prove an explicit formula for $\ffid$ we need two
lemmas. The second one will also be used in the next section.
\begin{lemma}
There exists a net $(x_l)_{l \in L}$ in $\gewn_0 \cap \gewn_0^*$ such that
$x_l$ is analytic w.r.t. $\sigma'$ for all $l$ and $\sigma'_z(x_l)
\recht 1$ strong* and bounded for all $z \in \complex$.
\end{lemma}
\begin{proof}
Because $\gewn_0 \cap \gewn_0^*$ is a strongly dense *-subalgebra of $\rondm$ we can take a
net $(a_k)_{k \in K}$ in $\gewn_0 \cap \gewn_0^*$ such that $a_k^*=a_k$,
$\|a_k\| \leq 1$ for all $k$ and $a_k \recht 1$ strongly. Define $q_k
\in \rondm$ by
\begin{equation*}
q_k = \frac{1}{\sqrt{\pi}} \int \exp(-t^2) \sigma'_t(a_k) \; dt.
\end{equation*}
Clearly $q_k$ is analytic w.r.t. $\sigma'$ and
\begin{equation*}
\sigma'_z(q_k) = \frac{1}{\sqrt{\pi}} \int \exp(-(t-z)^2)
\sigma'_t(a_k) \; dt.
\end{equation*}
Also $\sigma'_z(q_k) \recht 1$ strong* and bounded.

Define $L = \mathbb{N}_0 \times K \times \mathbb{N}_0$ with the product
order, and $x_{(n,k,m)} = e_n q_k e_m$. Then $x_l$ is analytic w.r.t.
$\sigma'$ for all $l$ and $\sigma'_z(x_l) \recht 1$ strong* and
bounded for all $z \in \complex$. Let $n,m \in \mathbb{N}_0$ and $k \in
K$. The operator $e_n q_k e_m \delta^{\et}$ is bounded, with closure
\begin{equation*}
e_n q_k (\delta^{\et} e_m) = e_n \frac{1}{\sqrt{\pi}} \int \exp(-t^2)
\sigma'_t(a_k)(\delta^{\et} e_m) \; dt.
\end{equation*}
For all $t \in \reeel$ the integrand of this expression equals
\begin{equation*}
\exp(-t^2)\delta^{it} \lambda^{\et it^2} \sigma_t(\overline{a_k
\delta^{\et}}(\lambda^{\et (it^2 - t)} \delta^{-it} \sigma_{-t}(e_m))).
\end{equation*}
This belongs to $\gewn$. When we apply $\Lambda$ on it we obtain
\begin{equation*}
\exp(-t^2) \delta^{it} \lambda^{\et it^2} \Delta^{it} J\lambda^{-\et it^2}
\delta^{it} \sigma_{-\iet -t}(e_m)J \Lambda(\overline{a_k
\delta^{\et}}).
\end{equation*}
As a function of $t$ this is weakly integrable. Because the mapping
$\Lambda$ is weak operator -- weak closed we can conclude that $e_n q_k
(\delta^{\et}e_m) \in \gewn$. This means that $e_n q_k e_m \in \gewn_0$.
Analogously we obtain that $e_n q_k e_m \in \gewn_0^*$.
\end{proof}
\begin{lemma} \label{lemma13}
If $a \in \gewn'$, then $a(\delta^z e_n)$ belongs to $\gewn$ for all
$z \in \complex$ and $n \in \mathbb{N}_0$. We have
$\Lambda(a(\delta^ze_n)) = \Lambda'(a(\delta^{z-\et}e_n))$.
\end{lemma}
\begin{proof}
Take a net $(x_l)_{l \in L}$ as in the previous lemma. Then
$a x_l (\delta^ze_n) \recht a (\delta^z e_n)$ strongly. Because $x_l \in
\gewn_0$ we have $a x_l (\delta^ze_n) \in \gewn$ and
\begin{align*}
\Lambda(a x_l (\delta^z e_n)) &= \Lambda'(a x_l (\delta^{z-\et} e_n))
\\ &= J'(\sigma'_{\iet}(x_l (\delta^{z-\et}e_n)))^*J'\Lambda'(a) \\
& \recht J'(\sigma'_{\iet} (\delta^{z-\et}e_n))^*J'
\Lambda'(a)=\Lambda'(a(\delta^{z-\et} e_n)).
\end{align*}
Because $\Lambda$ is weak operator -- weak closed, we conclude that
$a(\delta^z e_n) \in \gewn$ and $\Lambda(a(\delta^ze_n)) = \Lambda'(a(\delta^{z-\et}e_n))$.
\end{proof}
\begin{prop}
For all $x \in \rondm^+$ we have
\begin{equation*}
\ffid(x) = \lim_n \ffi ((\delta^{\et}e_n) x(\delta^{\et}e_n)).
\end{equation*}
\end{prop}
\begin{proof}
If $x \in \gewn'$ we have $x(\delta^{\et}e_n) \in \gewn$ for all $n$
and
\begin{align*}
\ffi((\delta^{\et}e_n) x^*x(\delta^{\et}e_n)) &= \|\Lambda(x
(\delta^{\et}e_n)) \|^2 = \| \Lambda'(xe_n) \|^2 \\
&= \|J' \sigma_{-\iet}(e_n)J' \Lambda'(x) \|^2 \recht
\|\Lambda'(x)\|^2 = \ffid(x^*x).
\end{align*}
This gives the proof for all $x \in {\gewm'}^+$. Now let $x \in
\rondm^+$ and $\ffid(x) = + \infty$. Suppose
$\ffi((\delta^{\et}e_n) x(\delta^{\et}e_n))$ does not converge to
$+\infty$. Then there exists a $M > 0$ and a subsequence $(e_{n_k})_k$
such that $\ffi ((\delta^{\et}e_{n_k}) x(\delta^{\et}e_{n_k})) \leq
M$ for all $k$. Thus $x^{\et} (\delta^{\et}e_{n_k}) \in \gewn$ for all
$k$, so $x^{\et}e_{n_k} \in \gewn'$ and $\ffid(e_{n_k} x e_{n_k}) \leq
M$ for all $k$. Because $e_{n_k} x e_{n_k} \recht x$ strong* and bounded, this
contradicts with $\ffid(x)= +\infty$ and the $\sigma$-weak lower
semicontinuity of $\ffid$.
\end{proof}
\section{The Connes cocycle $[D \ffid : D \ffi]$}
\begin{defin}
The expression $\lambda^{\et it^2} \delta^{it} \Delta^{it}$ defines a
strongly continuous one parameter group of unitaries on $\hilh$. So we
define the strictly positive operator $\rho$ such that $\rho^{it}$
equals this expression for all $t \in \reeel$.
\end{defin}
Recall the notation $S(z)$ from the end of the introduction.
\begin{lemma}
If $x \in \gewn \cap {\gewn'}^*$, then $\Lambda(x) \in
\mathcal{D}(\rho^{\et})$ and
\begin{equation*}
J \lambda^{-i/8}\rho^{\et} \Lambda(x) = \Lambda'(x^*).
\end{equation*}
\end{lemma}
\begin{proof}
Let $x \in \gewn \cap {\gewn'}^*$ and $n,m \in \mathbb{N}$. Then
$e_m x \in \gewn \cap \gewn^*$ because of lemma~\ref{lemma13}. We
can define a function from $S(\mbox{\small $\frac{1}{2}$})$ to
$\hilh$ mapping $\alpha$ to
\begin{equation*}
\lambda^{-\et i \alpha^2}
\delta^\alpha \sigma_{-i \alpha}(e_n) \Delta^\alpha \Lambda(e_mx).
\end{equation*}
This function is continuous on $S(\mbox{\small $\frac{1}{2}$})$ and analytic on its interior. It attains the
value
\begin{equation*}
\lambda^{\et it^2}\delta^{it}\sigma_t(e_n) \Delta^{it} \Lambda(e_mx)
= \rho^{it} \Lambda(e_n e_m x)
\end{equation*}
in $it$. So $\Lambda(e_n e_m x) \in \mathcal{D}(\rho^{\et})$ and
\begin{align*}
J\lambda^{-i/8} \rho^{\et} \Lambda(e_ne_mx) &= J
\sigma_{-\iet}(\delta^{\et}e_n) \Delta^{\et} \Lambda(e_mx) \\
&=J\sigma_{-\iet}(\delta^{\et}e_n)J \Lambda(x^*e_m) \\
&=\Lambda(x^*(\delta^{\et}e_n)e_m) = \Lambda'(x^*e_ne_m) \\
&=J' \sigma'_{-\iet}(e_n e_m)J' \Lambda'(x^*).
\end{align*}
Because $\rho^{\et}$ is closed we can conclude that $\Lambda(x) \in
\mathcal{D}(\rho^{\et})$ and $J \lambda^{-i/8}\rho^{\et} \Lambda(x) =
\Lambda'(x^*)$.
\end{proof}
\begin{prop} \label{Connes}
The Connes cocycle $[D\ffid : D \ffi]_t$ equals
$\lambda^{\et it^2}\delta^{it}$ for all $t \in \reeel$.
\end{prop}
\begin{proof}
Let $x \in \gewn^* \cap \gewn'$ and $y \in \gewn \cap {\gewn'}^*$.
Denote $u_t = \lambda^{\et it^2} \delta^{it}$. Define
$F(\alpha)=\langle \rho^\alpha \Lambda(y), \Lambda(x^*) \rangle$ when
$\alpha \in \complex$ and $0 \leq \text{Re}(\alpha) \leq \mbox{\small
$\frac{1}{2}$}$. Define $G(\alpha) = \langle \lambda^{i/8} J
\Lambda'(y^*), \rho^{\bar{\alpha}-1} \lambda^{i/8} J \Lambda'(x)
\rangle$ when $\alpha \in \complex$ and $\mbox{\small
$\frac{1}{2}$} \leq \text{Re}(\alpha) \leq 1$.
Because of the previous lemma $F$ and $G$ are both well defined,
continuous on their domain and analytic in the interior. For any $t
\in \reeel$ we have
\begin{align*}
F(it) &= \langle \lambda^{\et it^2} \delta^{it} \Delta^{it}
\Lambda(y),\Lambda(x^*) \rangle
=\ffi (x u_t \sigma_t(y)) \\
F(it + \mbox{\small
$\frac{1}{2}$}) &=\langle \rho^{it} \lambda^{i/8} J
\Lambda'(y^*),\Lambda(x^*) \rangle \\
G(it + \mbox{\small
$\frac{1}{2}$}) &= \langle \lambda^{i/8} J \Lambda'(y^*), \rho^{-it}
\Lambda(x^*) \rangle = F(it + \mbox{\small
$\frac{1}{2}$}) \\
G(it + 1) &= \langle J \Lambda'(y^*), \lambda^{\et it^2} \delta^{-it}
\Delta^{-it} J \Lambda'(x) \rangle
=\langle \Lambda'(x),J \lambda^{\et it^2}\delta^{it} J \Delta^{it}
\Lambda'(y^*) \rangle \\
&=\langle \lambda^{\et it^2} \delta^{it} \Lambda'(x), {\Delta'}^{it}
\Lambda'(y^*) \rangle = \ffid(\sigma'_t(y)u_t x).
\end{align*}
So we can glue together the functions $F$ and $G$ and define
$H(\alpha) = F(\alpha)$ when $\alpha$ belongs to the domain of $F$
and $H(\alpha) = G(\alpha)$ when $\alpha$ belongs to the domain of
$G$.
Then $H$ is continuous on $S(1)$ and analytic on its interior. We have
\begin{equation*}
H(it)=\ffi(x u_t \sigma_t(y)) \quad\text{and}\quad H(it+1) =
\ffid(\sigma'_t(y)u_t x)
\end{equation*}
for all $t \in \reeel$. Because it is easily verified that
\begin{align*}
u_{t+s} &=u_t \sigma_t(u_s) \\ u_{-t} &=\sigma_{-t}(u_t^*) \\
\sigma'_t(x) &=u_t \sigma_t(x) u_t^*
\end{align*}
for all $s,t \in \reeel$ and $x \in \rondm$, we conclude that $u_t =
[D \ffid : D \ffi]_t$ for all $t$.
\end{proof}
The previous proposition also implies that the operators $\lambda$
and $\delta$ are uniquely determined by $\ffi_\delta$. If we put $u_t
= [D \ffid : D \ffi]_t$ we have $\lambda^{it} = u_t^* u_1^* u_{t+1}$
and $\delta^{it} = u_t \lambda^{- \et i t^2}$ for all $t \in \reeel$,
which proves our claim.
\section{Three Radon-Nikodym theorems}
In this paragraph we denote by $(\sigma_t^\ffi)$ the modular
automorphism group of a n.s.f. weight $\ffi$ on a von Neumann algebra. We
denote by $\gewn_\ffi,\gewm_\ffi,\Lambda_\ffi,J_\ffi$ and
$\Delta_\ffi$ the same objects as defined in the introduction but we
add a subscript $\ffi$ for the sake of clarity.
\begin{prop} \label{Radon}
Let $\psi$ and $\ffi$ be two n.s.f. weights on a von Neumann algebra
$\rondm$. Let $\lambda$ and $\delta$ be two strongly commuting,
strictly positive operators affiliated with $\rondm$. Then the
following are equivalent
\begin{enumerate}
\item $[ D \psi : D \ffi]_t = \lambda^{\et it^2} \delta^{it}
\quad\text{for all}\quad t \in \reeel.$ \label{one}
\item $\sigma_t^\ffi(\delta^{is}) = \lambda^{ist}\delta^{is}$ for all $s,t
\in \reeel \quad\text{and}\quad \psi = \ffid$. \label{two}
\end{enumerate}
\end{prop}
\begin{proof}
The implication \ref{two}) $\Rightarrow$ \ref{one}) follows from the proposition~\ref{Connes}.

To prove \ref{one}) $\Rightarrow$ \ref{two}) denote $u_t = [D \psi : D \ffi]_t$. Let
$s,t \in \reeel$. Then
\begin{equation*}
\lambda^{\et it^2} \lambda^{\et is^2} \lambda^{ist} \delta^{it}
\delta^{is} = u_{t+s} = u_t \sigma_t^\ffi(u_s) = \lambda^{\et it^2}\delta^{it}
\sigma_t^\ffi (\lambda^{\et is^2} \delta^{is}).
\end{equation*}
This implies that
\begin{equation} \label{eerste}
\lambda^{\et is^2}\lambda^{ist} \delta^{is} =
\sigma_t^\ffi(\lambda^{\et is^2}\delta^{is}) \quad\text{for all} \quad s,t
\in \reeel.
\end{equation}
It follows that for all $r,s,t \in \reeel$
\begin{equation} \label{tweede}
\sigma_r^\ffi(\lambda^{\et is^2 + ist}) \sigma_r^\ffi(\delta^{is}) =
\sigma_{r+t}^\ffi(\lambda^{\et is^2} \delta^{is}) =
\lambda^{\et is^2}\lambda^{is(r+t)} \delta^{is}.
\end{equation}
But equation~(\ref{eerste}) implies that $\sigma_r^\ffi(\delta^{is}) =
\sigma_r^\ffi(\lambda^{-\et is^2}) \lambda^{\et is^2}
\lambda^{isr}\delta^{is}$, so by equation~(\ref{tweede}) we get
\begin{equation*}
\sigma_r^\ffi(\lambda^{ist}) \lambda^{\et is^2} \lambda^{isr} \delta^{is}
= \lambda^{\et is^2} \lambda^{is(r+t)} \delta^{is}.
\end{equation*}
This gives us $\sigma_r^\ffi(\lambda^{ist}) = \lambda^{ist}$ for all $r,s,t
\in \reeel$. Then equation~(\ref{eerste}) implies that
$\sigma_t^\ffi(\delta^{is}) = \lambda^{ist}\delta^{is}$ for all $s,t \in
\reeel$. So we can construct the weight $\ffid$, such that $[D \ffid
: D \ffi]_t = \lambda^{\et it^2}\delta^{it}$. But then $\ffid = \psi$.
\end{proof}
We will now consider the more specific case in which $\lambda$ is
affiliated to the centre of $\rondm$. We will prove that we obtain
exactly all the weights whose automorphism group commutes with that
of $\ffi$.
\begin{prop}
Let $\ffi$ and $\psi$ be two n.s.f. weights on a von Neumann algebra
$\rondm$. Then the following are equivalent.
\begin{enumerate}
\item The modular automorphism groups $\sigma^\psi$ and $\sigma^\ffi$
commute. \label{rn1}
\item There exist a strictly positive operator $\delta$ affiliated
with $\rondm$ and a strictly positive operator $\lambda$ affiliated
with the centre of $\rondm$ such that
$\sigma_s^\ffi(\delta^{it})=\lambda^{ist} \delta^{it}$ for all $s,t
\in \reeel$ and such that $\psi = \ffi_\delta$. \label{rn2}
\item There exist a strictly positive operator $\delta$ affiliated
with $\rondm$ and a strictly positive operator $\lambda$ affiliated
with the centre of $\rondm$ such that $[D \psi : D \ffi]_t =
\lambda^{\et it^2} \delta^{it}$ for all $t \in \reeel$. \label{rn3}
\end{enumerate}
\end{prop}
\begin{proof}
The equivalence of~\ref{rn2}) and \ref{rn3}) follows from
proposition~\ref{Radon}. The implication
\ref{nr2})~$\Rightarrow$~\ref{nr1}) follows from corollary~\ref{Cor} by
a direct computation. We will prove the implication
\ref{nr1})~$\Rightarrow$~\ref{nr3}). Denote $u_t = [D\psi : D
\ffi]_t$ for all $t \in \reeel$ and denote by $\rondz$ the centre of
$\rondm$. For all $x \in \rondm$ and $s,t \in \reeel$ we have
\begin{equation*}
\sigma_t^\psi(\sigma_s^\ffi(x)) = u_t \sigma_{t+s}^\ffi(x) u_t^*
\quad\text{and}\quad \sigma_s^\ffi(\sigma_t^\psi(x)) =
\sigma_s^\ffi(u_t) \sigma_{t+s}^\ffi(x) \sigma_s^\ffi(u_t^*).
\end{equation*}
Thus we can conclude that $u_t^* \sigma_s^\ffi(u_t) \in \rondz$ for
all $s, t \in \reeel$. But then $u_t^* u_s^* u_{s+t} \in \rondz$ for
all $s, t \in \reeel$. Because $\sigma^\ffi$ acts trivially on
$\rondz$ we get
\begin{equation*}
u_t^* \sigma_s^\ffi(u_t) = \sigma_{-t}^\ffi(u_t^* \sigma_s^\ffi(u_t))
=u_{-t} \sigma_{s-t}^\ffi(u_t) = u_{-t}u_{s-t}^*u_s \in \rondz.
\end{equation*}
We can conclude that $u_t u_{s+t}^* u_s \in \rondz$ for all $s,t \in
\reeel$. Then we define for $(s,t) \in \reeel^2$, $w(s,t)=u_t^* u_s^*
u_{s+t}$. The function $w$ is strong* continuous from $\reeel^2$ to
the unitaries of $\rondz$. Let $s,s',t \in \reeel$. Because of the
previous remarks we can make the following calculation.
\begin{align*}
w(s+s',t) &= u_t^* \: u_{s+s'}^* \: u_{s+s'+t} \\
&= u_t^* \: \sigma_{s'}^\ffi(u_s^*) \: u_{s'}^* \: u_{s'+t} \: u_s \: (u_s^*
\: \sigma_{s'+t}^\ffi(u_s)) \\
&=u_t^* \: u_s^* \: \sigma_{s'}^\ffi(\sigma_t^\ffi(u_s) \: u_s^*) \: u_{s'}^*
\: u_{s'+t} \: u_s \\
&=u_t^* \: u_s^* \: \sigma_{s'}^\ffi(u_t^* \: u_{s+t} \: u_s^*) \: u_{s'}^*
\: u_{s'+t} \: u_s \\
&=u_t^* \: u_s^* \: (u_{s'}^* \: u_{s'+t} \: u_t^*) \: u_{s+t} \: (u_t^* \: u_t) \\
&=u_t^* \: u_s^* \: u_{s+t} \: u_t^* \: u_{s'}^* \: u_{s'+t} = w(s,t) \: w(s',t).
\end{align*}
Next let $s,t,t' \in \reeel$. We have
\begin{align*}
w(s,t+t') &=u_{t+t'}^* \: u_s^* \: u_{s+t+t'} \\
&=(\sigma_{t'}^\ffi(u_t^*) \: u_t) \: u_t^* \: u_{t'}^* \: u_s^* \: u_{t'}
\: \sigma_{t'}^\ffi(u_{s+t}) \\
&=u_t^* \: u_{t'}^* \: u_s^* \: u_{t'} \: \sigma_{t'}^\ffi(u_{s+t} \: u_t^*) \: u_t \\
&=u_t^* \: u_{t'}^* \: u_s^* \: u_{t'} \: \sigma_{t'}^\ffi(u_s) \: \sigma_{t'}^\ffi
(u_s^* \: u_{s+t} \: u_t^*) \: u_t \\
&=u_t^* \: u_{t'}^* \: u_s^* \: u_{t'} \: \sigma_{t'}^\ffi(u_s) \: u_s^* \: u_{s+t} \\
&=u_t^* \: (u_{t'}^* \: u_s^* \: u_{s+t'}) \: u_s^* \: u_{s+t} \\
&=u_t^* \: u_s^* \: u_{s+t} \: u_{t'}^* \: u_s^* \: u_{s+t'} = w(s,t) \: w(s,t').
\end{align*}
For each $t \in \reeel$ we can now take a strictly positive operator
$\lambda_t$ affiliated with $\rondz$ such that $\lambda_t^{is}
=w(s,t)$ for all $s,t \in \reeel$. Let $t,t' \in \reeel$. Because
$\lambda_t$ and $\lambda_{t'}$ are strongly commuting we can write
\begin{equation*}
(\lambda_t \hat{\cdot} \lambda_{t'})^{is} = \lambda_t^{is}
\lambda_{t'}^{is} = w(s,t)w(s,t') = \lambda_{t+t'}^{is}
\end{equation*}
for all $s \in \reeel$, where $\lambda_t \hat{\cdot} \lambda_{t'}$
denotes the closure of $\lambda_t \lambda_{t'}$. It follows that
$\lambda_{t+t'} = \lambda_t \hat{\cdot} \lambda_{t'}$. Put
$\lambda=\lambda_1$. It follows from functional calculus that
$\lambda^q = \lambda_q$ for all $q \in \mathbb{Q}$. Then we have
\begin{equation*}
\lambda^{isq} = (\lambda^q)^{is} = \lambda_q^{is} = w(s,q)
\end{equation*}
for all $s \in \reeel$ and $q \in \mathbb{Q}$. Because of strong*
continuity we have $\lambda^{ist} = w(s,t)$ and thus $u_{s+t} =
\lambda^{ist}u_t u_s$ for all $s,t \in \reeel$. Now we can easily
verify that $v_t = \lambda^{-\et it^2} u_t$ defines a strong*
continuous one-parameter group of unitaries in $\rondm$. So we can
take a strictly positive operator $\delta$ affiliated with $\rondm$
such that $[D \psi : D \ffi]_t = u_t = \lambda^{\et i t^2}
\delta^{it}$ for all $t \in \reeel$. This gives us~\ref{rn3}).
\end{proof}
Now we will look at the even more specific case $\lambda \in
\reeel_0^+$.
So the following proposition becomes meaningful.
\begin{prop}
Let $\ffi$ be a n.s.f. weight on a von Neumann algebra $\rondm$. Let
$\delta$ be a strictly positive operator affiliated with $\rondm$ and
$\lambda \in \reeel_0^+$ such that $\sigma_t^\ffi(\delta^{is}) =
\lambda^{ist} \delta^{is}$ for all $s,t \in \reeel$. Then we have
\begin{equation*}
\ffid \circ \sigma_t^\ffi = \lambda^{-t} \ffid \qquad\text{and}\qquad \ffi
\circ \sigma_t^{\ffid} = \lambda^t \ffi \qquad\text{for all}\quad t
\in \reeel.
\end{equation*}
\end{prop}
\begin{proof}
Let $a \in \gewn_{\ffid}$ and $t \in \reeel$. Then $\sigma_t^\ffi(a) = \delta^{-it} \sigma_t^{\ffid}(a)
\delta^{it}$. This belongs to $\gewn_{\ffid}$ because $\delta^{it}$ is
analytic w.r.t. $\sigma^{\ffid}$, and we have
\begin{equation*}
\Lambda_{\ffid}(\sigma_t^\ffi(a)) = \delta^{-it} J_{\ffid} \lambda^{-\et t} \delta^{-it}
J_{\ffid} \Delta_{\ffid}^{it} \Lambda_{\ffid}(a).
\end{equation*}
So we get
\begin{equation*}
\ffid(\sigma_t^\ffi(a^*a)) = \lambda^{-t} \ffid(a^*a) \quad\text{for
all}\quad t \in \reeel.
\end{equation*}
Now, the conclusion follows easily. The second statement is
proved analogously.
\end{proof}
After stating a lemma we will prove our third Radon-Nikodym theorem.
\begin{lemma}
Let $\ffi$ be a n.s.f. weight on a von Neumann algebra $\rondm$ and $a \in \rondm$. If
$\gewn_\ffi a \subset \gewn_\ffi$, $\gewn_\ffi a^* \subset \gewn_\ffi$
and if there exists a $\lambda \in \reeel_0^+$ such that
$\ffi(ax)=\lambda \ffi(xa)$ for all $x \in \gewm_\varphi$, then
$\sigma_t^\varphi(a)=\lambda^{it} a$ for all $t \in \reeel$.
\end{lemma}
\begin{proof}
The proof of Result~6.29 in \cite{kustermans1} can be taken over
literally. Also a slight adaptation of the proof of theorem~3.6 in
\cite{pedersen} yields the result.
\end{proof}
\begin{prop}
Let $\psi$ and $\ffi$ be two n.s.f. weights on a von Neumann algebra
$\rondm$. Let $\lambda \in \reeel_0^+$. The following statements are
equivalent.
\begin{enumerate}
\item For all $t \in \reeel$ we have $\ffi \circ \sigma_t^\psi =
\lambda^t \ffi$. \label{nr1}
\item For all $t \in \reeel$ we have $\psi \circ \sigma_t^\ffi =
\lambda^{-t} \psi$. \label{nr2}
\item There exists a strictly positive operator $\delta$ affiliated
with $\rondm$ such that $\sigma_t^\ffi(\delta^{is}) = \lambda^{ist}
\delta^{is}$ for all $s,t \in \reeel$ and such that $\psi = \ffid$.
\label{nr3}
\item There exists a strictly positive operator $\delta$ affiliated
with $\rondm$ such that $[D\psi : D \ffi]_t = \lambda^{\et it^2}
\delta^{it}$ for all $t \in \reeel$. \label{nr4}
\end{enumerate}
\end{prop}
\begin{proof}
We have already proven the equivalence of \ref{nr3}) and \ref{nr4}) and
the implications \ref{nr3})~$\Rightarrow$~\ref{nr2}) and
\ref{nr3})~$\Rightarrow$~\ref{nr1}). Suppose now that \ref{nr1}) is valid.
Put $u_t = [D\psi : D \ffi]_t$ and let $x \in \rondm^+$. Then we have
\begin{align*}
\ffi(u_t^*xu_t) &= \ffi(\sigma_{-t}^\ffi(u_t^*) \sigma_{-t}^\ffi(x)
\sigma_{-t}^\ffi(u_t)) \\&=\ffi(u_{-t} \sigma_{-t}^\ffi(x) u_{-t}^*)
= \ffi(\sigma_{-t}^\psi(x)) = \lambda^{-t} \ffi(x).
\end{align*}
So we have $\gewn_\ffi u_t \subset \gewn_\ffi$ for all $t \in
\reeel$, and thus $\gewn_\ffi u_t^* = \sigma_t^\ffi(\gewn_\ffi
u_{-t}) \subset \gewn_\ffi$ for all $t \in \reeel$. Then we get for
every $x \in \gewm_\ffi$ that
\begin{equation*}
\ffi(xu_t) = \ffi(u_t^*u_t x u_t) = \lambda^{-t} \ffi(u_t x).
\end{equation*}
From the previous lemma we can conclude that $\sigma_s^\ffi(u_t) =
\lambda^{ist} u_t$ for all $s,t \in \reeel$. Put $v_t =
\lambda^{-\et it^2} u_t \in \rondm$. Then we have that $t \mapsto v_t$
is a strongly continuous one-parameter group of unitaries. Define
$\delta$ such that $\delta^{it} = v_t$ for all $t \in \reeel$. So
$\delta$ is affiliated with $\rondm$ and $[D\psi : D \ffi]_t =
\lambda^{\et it^2} \delta^{it}$ and this gives us \ref{nr4}). Finally
suppose \ref{nr2}) is valid. From the proven implication
\ref{nr1})~$\Rightarrow$~\ref{nr4}) we get the existence of a strictly
positive operator $\delta$ affiliated with $\rondm$ such that $[D\ffi
: D \psi]_t = \lambda^{-\et it^2} \delta^{it}$ for all $t \in \reeel$.
Changing $\delta$ to $\delta^{-1}$ we get $[D\psi : D\ffi]_t =
\lambda^{\et it^2} \delta^{it}$ for all $t \in \reeel$. This gives us
again \ref{nr4}).
\end{proof}
\newcommand{\Trace}{\operatorname{Tr}}
We conclude this paper by giving an example that shows all situations
can really occur : we can have $\sigma_t^\ffi(\delta^{is}) =
\lambda^{ist} \delta^{is}$ with $\lambda$ and $\delta$ strongly
commuting but $\lambda$ not central, with $\lambda$ central but not
scalar, and with $\lambda$ scalar. Indeed, define $\rondm_1 =
B(L^2(\reeel))$ and define the selfadjoint operators $P$ and $Q$ on
the obvious domains by
\begin{equation*}
(P \xi)(\gamma) = \gamma \xi(\gamma) \quad \text{and} \quad (Q
\xi)(\gamma) = -i \xi'(\gamma).
\end{equation*}
Put $H = \exp(P)$ and $K_1 = \exp(Q)$ and denote by $\Trace$ the
canonical trace on $\rondm_1$. Remark that $\Trace$ has a trivial modular
automorphism group such that we can define $\ffi_1 = \Trace_H$ as in definition~\ref{def15}. An easy
calculation yields that $\sigma_t^{\ffi_1}(K_1^{is}) = H^{it} K_1^{is} H^{-it}= e^{-its}
K_1^{is}$, where $e$ denotes the well known real number $e$. This gives
an example of our third case. Define
$\rondm_2$ as the von Neumann algebra of two by two matrices over
$\rondm_1$ and $\ffi_2$ as the balanced weight
$\theta(\ffi_1,\ffi_1)$ (see \cite{stratila2}). Define $K_2 = \left( \begin{smallmatrix}
K_1 & 0 \\ 0 & K_1^{-1} \end{smallmatrix} \right)$
We easily have $\sigma_t^{\ffi_2}(K_2^{is}) =
\left( \begin{smallmatrix}
e^{-1} & 0 \\ 0 & e \end{smallmatrix} \right)^{its} K_2^{is}$,
which gives an example of our first
case because $\rondm_2$ is a factor. Define $\rondm_3$ as the
diagonal matrices in $\rondm_2$. We can restrict $\ffi_2$ to
$\rondm_3$ and keep $K_2$. We have the same formula as above, and
this way an example of our second case, $\left( \begin{smallmatrix}
e^{-1} & 0 \\ 0 & e \end{smallmatrix} \right)$ being
central now.
\clearpage
\providecommand{\bysame}{\leavevmode\hbox to3em{\hrulefill}\thinspace}


\begin{thebibliography}{1}
  \addcontentsline{toc}{chapter}{\numberline{}Bibliografie}

\bibitem{kustermans1}
J.~Kustermans, \emph{{KMS-}weights on {$C^*$}-algebras}, preprint {K.U.Leuven},
  {$\#$ Funct-an/9704008}, 1997.

\bibitem{kustvaes}
J.~Kustermans and S.~Vaes, \emph{Reduced {$C^*$}-algebraic quantum groups}, in
  preparation, 1998.

\bibitem{kustermans2}
J.~Kustermans and A.~Van~Daele, \emph{{$C^*$-algebraic} quantum groups arising
  from algebraic quantum groups}, Int. J. Math. \textbf{8} (1997), no.~8,
  1067--1139.

\bibitem{masuda}
T.~Masuda and Y.~Nakagami, \emph{A von {N}eumann algebra framework for the
  duality of the quantum groups}, Publications of the RIMS Kyoto University
  \textbf{30} (1994), 799--850.

\bibitem{pedersen}
G.~K. Pedersen and M.~Takesaki, \emph{The {R}adon-{N}ikodym theorem for von
  {N}eumann algebras}, Acta Math. \textbf{130} (1973), 53--87.

\bibitem{stratila2}
S.~Stratila, \emph{Modular theory in operator algebras}, Abacus Press,
  Turnbridge Wells, England, 1981.

\bibitem{stratila}
S.~Stratila and L.~Zsido, \emph{Lectures on von {N}eumann algebras}, Abacus
  Press, Turnbridge Wells, England, 1979.

\end{thebibliography}
\end{document}